\documentclass{amsart}      
\usepackage{amssymb, amsmath}
\newtoks\prt 
 
\newtheorem{thm}{Theorem}
 
\newtheorem{lemma}[thm]{Lemma}

\newtheorem{example}[thm]{Example}

\newtheorem*{question}{Question} 
\theoremstyle{definition} 
 
\newtheorem{remark}[thm]{Remark}

\def\eqn#1$$#2$${\begin{equation}\label#1#2\end{equation}}

\def\ep{\varepsilon} 
 
\def\en{\mathbb N} 
\def\er{\mathbb R}

\def\wh{\widehat} 
\def \reg {\partial _{\kern1pt\text{reg}}}

\def\dd{\operatorname{d}}
\def\dh{\widehat{\operatorname{d}}}
\def\clu#1#2{\operatorname{clust}_{#1^{**}}(#2)}
\def\wde#1{\widetilde{\delta}\left(#1\right)}
\def\de#1#2{\delta\left(#2\right)}

\newcommand{\norm}[1]{\left\|#1\right\|}

\newcommand{\betr}[1]{| #1  |}
\renewcommand{\Re}{\operatorname{Re}}

\begin{document}

\title{On quantification of weak sequential completeness}
\author{O.F.K. Kalenda, H. Pfitzner and J. Spurn\'y}

\address{Department of Mathematical Analysis \\
Faculty of Mathematics and Physic\\ Charles University\\
Sokolovsk\'{a} 83, 186 \ 75\\Praha 8, Czech Republic}
\email{kalenda@karlin.mff.cuni.cz}
\email{spurny@karlin.mff.cuni.cz}

\address{Universit\'e d'Orl\'eans\\
BP 6759\\
F-45067 Orl\'eans Cedex 2\\France}
\email{hermann.pfitzner@univ-orleans.fr}


\subjclass[2010]{46B20}
\keywords{weakly sequentially complete Banach space; $L$-embedded Banach space; quantitative versions of weak sequential completeness}

\thanks{The first and third authors were supported in part by the grant
GAAV IAA 100190901  and in part by the Research Project MSM~0021620839 from the Czech Ministry of Education.}

\begin{abstract} 
We consider several quantities related to weak sequential completeness of a Banach space
and prove some of their properties in general and in $L$-embedded Banach spaces,
improving in particular an inequality of G.~Godefroy, N.~Kalton and D.~Li.
We show some examples witnessing natural limits of our positive results,
in particular, we construct a separable Banach space $X$ with the Schur property that cannot be renormed
to have a certain quantitative form of weak sequential completeness,
thus providing a partial answer to a question of G.~Godefroy.
\end{abstract}
\maketitle


\section{Introduction and statement of the results}

If $X$ is a Banach space, we recall that it is \emph{weakly sequentially complete} if any weakly Cauchy sequence in $X$ is weakly convergent.
In the present paper we investigate quantitative versions of this property.
To this end we use several quantities related to a given bounded sequence $(x_k)$ in $X$.

Let $\clu{X}{x_k}$ denote the set of all weak* cluster points of $(x_k)$ in $X^{**}$. By $\de{X}{x_k}$
we will denote the diameter of $\clu{X}{x_k}$ (see also \eqref{gl} below).
Further, if  $A$, $B$ are  nonempty subsets of a Banach space $X$, then $\operatorname{d}(A,B)$ denotes the
usual distance between $A$ and $B$  and the Hausdorff non-symmetrized distance from $A$ to $B$ is defined by
\[
   \dh(A,B)= \sup\{\operatorname{d}(a,B): a\in A\}.
\]

Note that a space $X$ is weakly sequentially complete if for each bounded sequence $(x_k)$ in $X$ satisfying $\de{X}{x_k}=0$
(this just means that the sequence is weakly Cauchy) we have
$\dh(\clu{X}{x_k},X)=0$ (i.e., all the weak* cluster points are contained in $X$, which for a weakly Cauchy sequence means
that it is weakly convergent). It is thus natural to ask which Banach spaces satisfy a quantitative version of weak sequential completeness, i.e., the inequality
\begin{equation}\label{quantitative}\dh(\clu{X}{x_k},X)\leq C\cdot \de{X}{x_k}\end{equation}
for all bounded sequences $(x_k)$ in $X$ and for some $C>0$.
The starting point of our investigation was the following remark made by G.~Godefroy in \cite[p.\,829]{gode-handbook}:

\emph{If $X$ is complemented in $X^{**}$ by a projection $P$ satisfying 
\begin{equation}
\label{gl_L-proj}
\|x^{**}\|=\|Px^{**}\|+\|x^{**}-Px^{**}\|, \quad x^{**}\in X^{**},
\end{equation}
then $X$ is weakly sequentially complete and
\begin{equation}
\label{est-L1}
\dh(\clu{X}{x_k}, X)\leq \de{X}{x_k}
\end{equation}
for any sequence $(x_k)$ in $X$.
}

It can be easily seen that
\begin{equation}
\label{gl}
\begin{aligned}
\de{X}{x_k}&=\sup_{x^*\in B_{X^*}} \left(\limsup_{k\to\infty} x^*(x_k)-\liminf_{k\to\infty} x^*(x_k)\right)\\&=\sup_{x^*\in B_{X^*}}\lim_{n\to\infty} \sup\{|x^*(x_l)-x^*(x_j)|: l,j\geq n\}
\end{aligned}
\end{equation}

The first formula of \eqref{gl} is used in \cite[Section~2.1]{behrends}, the second one in \cite[p.\,829]{gode-handbook}.

Banach spaces satisfying assumption \eqref{gl_L-proj} above are called \emph{$L$-embedded}, see \cite[Section~III.1]{hawewe}.
The proof of \eqref{est-L1} can be found in \cite[Lemma IV.7]{GoKaLi2}.

By what has been said above, inequality~\eqref{est-L1} is a quantitative form of weak sequential completeness.
In \cite[p.\,829]{gode-handbook} G.\ Godefroy mentions that it is not clear
which weakly sequentially complete spaces can be renormed to have such a quantitative form of weak sequential completeness.

The aim of our paper is twofold.
On the one hand we show that the answer to G.\ Godefroy's question cannot be positive for all weakly sequentially complete Banach spaces,
more precisely we construct a weakly sequentially complete space that cannot be renormed in such a way that \eqref{est-L1} holds, see
Example \ref{exa2} below.
On the other hand we put inequality~\eqref{est-L1} into context
by studying some modifications and possible converses, see the following theorem.
In particular, we slightly improve inequality~\eqref{est-L1} - see \eqref{est-L1bis} in the theorem -
but such that now the additional factor $2$ is optimal.

We will use one more quantity (cf.\ \cite{pfitzner-inve} but
appearing implicitely in \cite{behrends}) which in some situations is more natural than the quantity $\delta_X$, namely
\[
\wde{x_k}=\inf\left\{\de{X}{x_{k_j}}\colon (x_{k_j})\text{ is a subsequence of }(x_k)\right\}.
\]

\begin{thm}\label{positive} Let $X$ be a Banach space and $(x_k)$ be a bounded sequence in $X$. Then
\begin{equation}
\label{ceje}
\wde{x_k}\leq 2\dh(\clu{X}{x_k}, X).
\end{equation}
If the space $X$ is $L$-embedded, then also the following inequalities hold:
\begin{eqnarray}
\label{est-L1bis}
2\dh(\clu{X}{x_k}, X)\leq \de{X}{x_k}, \\
\label{est-L1tres}
2\operatorname{d}(\clu{X}{x_k},X)\leq \wde{x_k}.
\end{eqnarray}
\end{thm}
Since we have trivially that $\widetilde\delta\leq\delta_X$ and $\operatorname{d}\leq\dh$ it is natural to ask whether one of these
quantities can be replaced by a sharper one in the inequalities of the theorem.
The following remark and Example \ref{exa1} show that this cannot be done in any of the inequalities \eqref{ceje}-\eqref{est-L1tres}.
\begin{remark}\label{Rem}
(a) In \eqref{est-L1bis}, $\delta_X$ cannot be replaced by $\widetilde\delta$
and in \eqref{est-L1tres} $\operatorname{d}$ cannot be replaced by $\dh$.
This is witnessed by the sequence $(x_k)$ in $X=\ell_1$ such that $x_{2k-1}=0$ and $x_{2k}=e_k$ for all $k\in\en$. Then
$\operatorname{d}(\clu{X}{x_k},X)=\wde{x_k}=0$,  $\dh(\clu{X}{x_k}, X)=1$ and $\de{X}{x_k}=2$.\\
(b) Inequality \eqref{ceje} is a kind of converse of \eqref{est-L1} and holds in all Banach spaces. We note that
$\widetilde\delta$ cannot be replaced by $\delta_X$ in \eqref{ceje}, in other words,
inequality \eqref{est-L1} cannot be reversed as it is, neither in $L$-embedded spaces.
Indeed, let  $X=\ell_1$. We consider the elements $x_k=\frac1k e_k$ and $y_k=e_1+\frac1k e_k$, $k\in\en$.
Let $(z_k)$ be the sequence $x_1,y_1,x_2,y_2,\dots$. Then $\dh(\clu{\ell_1}{z_k}, \ell_1)=0$
because all weak$^*$ cluster points of $(z_k)$ are contained in $\ell_1$, but 
\[
\de{\ell_1}{z_k}\geq \limsup_{k\to\infty} e_1(z_k)-\liminf_{k\to\infty} e_1(z_k)=1.
\]
(c) We further remark that in all inequalities in Theorem~\ref{positive} the factor $2$ is optimal,
as witnessed by the sequence $(e_k)$ in $X=\ell_1$. Indeed, then
$$
\dh(\clu{X}{e_k}, X)=\operatorname{d}(\clu{X}{e_k},X)=1
\mbox{\qquad and\qquad}\wde{e_k}=\de{X}{e_k}=2.
$$
\end{remark}\medskip
It is also natural to ask whether $\dh$ can be replaced by $\operatorname{d}$ in the inequality~\eqref{ceje},
i.e., whether the inequality~\eqref{est-L1tres} can be reversed (at least for $L$-embedded spaces).
This is not the case by the following example. 

\begin{example}\label{exa1}
There is an $L$-embedded space $X$ and a bounded sequence $(x_k)$ in $X$ such that $\wde{x_k}=2$ and $\dd(\clu{X}{x_k},X)=0$.
\end{example}

The negative partial answer to the mentioned question of G.\ Godefroy is given by the following example. In fact, we obtain a slightly stronger result. Not only there is a weakly sequentially complete Banach space not satisfying \eqref{quantitative} for all bounded sequences and some $C>0$, but we get even a weakly sequentially complete space not satisfying a weaker form of \eqref{quantitative} -- with $\dd$ in place of $\dh$.

\begin{example}\label{exa2}
There exists a separable Banach space $X$ with the Schur property - in particular, $X$ is weakly sequentially complete - which is $1$-complemented in its bidual, 
such that there is no constant $C>0$ satisfying 
\[
\dd(\clu{X}{x_k}, X)\leq C\cdot\de{X}{x_k}
\]
for every bounded sequence $(x_k)$ in $X$.
\end{example}

\section{Proof of Theorem~\ref{positive}}

The proof relies on two simple properties of $\ell_1$-sequences which are formulated in the following lemma.

\begin{lemma}\label{lemma-l1} 
Let $X$ be a Banach space and $(x_n)$ be a bounded sequence in $X$. Suppose that $c>0$ is such that
$$\left\|\sum_{j=1}^n \alpha_j x_j\right\|\ge c\sum_{j=1}^n |\alpha_j|$$
whenever $n\in\en$ and $\alpha_1,\dots,\alpha_n$ are real numbers. Then
\begin{itemize}
	\item[(i)] $\de{X}{x_n}\ge 2c$,
	\item[(ii)] $\operatorname{d}(\clu{X}{x_k},X)\ge c$.
\end{itemize}
\end{lemma}

{\em Proof.} 
(i) It is clear that the sequence $(x_n)$ is linearly independent. Hence there is
a unique linear functional defined on its linear span whose value is $c$ at $x_{2k-1}$ and $-c$ at $x_{2k}$ for each $k\in\en$. By the assumption, the norm of this functional is at most $1$. Let $x^*\in B_{X^*}$ be its Hahn-Banach extension. Then $x^*$ witnesses that $\de{X}{x_n}\ge 2c$.

(ii) Let $x^{**}$ be any  weak* cluster point of the sequence $(x_n)$ in $X^{**}$ and $x\in X$ be arbitrary.
It follows from \cite[Proposition 4.2]{KnOd} that there is an index $m\in\en$ such that
\[
\left\|\sum_{j=m}^\infty \alpha_j(x_j-x)\right\|\geq c\sum_{j=m}^\infty |\alpha_j|
\]
for every sequence $(\alpha_j)_{j=m}^\infty$ 
with finitely many nonzero elements. In particular, it follows that the vectors $x_j-x$, $j\ge m$, are linearly independent. So, there is a unique linear functional on their linear span whose value at each $x_j-x$ is equal to $c$.
By the above inequality, the norm of this functional is at most one. Let $x^*\in X^*$ be its Hahn-Banach extension. Then we have
$$\|x^{**}-x\|\ge(x^{**}-x)(x^*)\ge\liminf_{j\to\infty}  x^*(x_j-x)=c.$$
This completes the proof of the lemma.\hfill$\Box$

\medskip

Now we are ready to prove Theorem~\ref{positive}:

We start by proving \eqref{ceje}: 
Let $(x_k)$ be a bounded sequence in $X$. We assume that $\wde{x_k}>0$ because otherwise~\eqref{ceje} holds trivially. 
Let $c\in(0,\wde{x_k})$ be arbitrary. The key ingredient is provided by a result of E.~Behrends (see \cite[Theorem 3.2]{behrends}) that yields a subsequence $(x_{n_k})$ such that
$$\norm{\sum_{i=1}^k \alpha_i x_{n_i}}\ge \frac c2 \sum_{i=1}^k|\alpha_i|$$
whenever $k\in\en$ and $\alpha_1,\dots,\alpha_k\in\er$. 
By Lemma~\ref{lemma-l1}(ii) we get
$\dd(\clu{X}{x_{n_k}},X)\geq \frac c2$, hence
$\dh(\clu{X}{x_{k}},X)\geq \frac c2$. As $c\in (0,\wde{x_k})$ is arbitrary,
 \eqref{ceje} follows.

We continue by proving \eqref{est-L1bis}: 
We set $c=\dh(\clu{X}{x_k}, X)$ and assume that $c>0$ because otherwise \eqref{est-L1bis} holds trivially.
Let $\ep\in(0,c)$ be arbitrary and let $x^{**}$ be a weak* cluster point of the sequence $(x_k)$ in $X^{**}$
such that $\operatorname{d}(x^{**},X)>c-\frac\ep2$.
Set $x=Px^{**}$ and $x_s=x^{**}-x$ where $P$ denotes the projection on $X$ as in \eqref{gl_L-proj}.
Then $\operatorname{d}(x^{**},X)=\norm{x_s}$. We claim that there is a subsequence $(x_{k_n})$ such that
\begin{equation}\label{gl7}
\norm{\sum_{i=1}^n\alpha_i (x_{k_i}-x)}\geq(c-(1-2^{-n})\ep)\sum_{i=1}^n\betr{\alpha_i}
\end{equation}
for all $n\in\en$ and all $(\alpha_i)_{i=1}^n$ in $\er^n$.
\smallskip
This will be proved by G.\ Godefroy's \textquoteleft ace of $\Diamond$ argument\textquoteright{} \cite[p.\ 170]{hawewe},
cf.\ the proof of \cite[Proposition IV.2.5]{hawewe}.
Since $x_s$ is a weak* cluster point of the sequence $(x_k-x)$, there is $k_1$ such that 
$\norm{x_{k_1}-x}>c-\frac\ep2$ which settles the first induction step.

Suppose we have constructed $x_{k_1}, \ldots, x_{k_n}$.
Let $(\alpha^l)_{l=1}^L$ be a finite sequence of elements of the unit sphere of 
$\ell_1^{n+1}$ such that $\alpha_{n+1}^l\neq0$ for all $l\in\{1,\dots,L\}$
and such that for each $\alpha$ in the unit sphere of $\ell_1^{n+1}$ there is an element $\alpha^l$ such that 
$$
\norm{\alpha-\alpha^l}_{\ell_1^{n+1}}<\frac{\ep}{2^{n+2}\sup_k\norm{x_k}}.
$$

Let $l\in\{1,\dots,L\}$ be arbitrary. Then
$\sum_{i=1}^n\alpha_i^l (x_{k_i}-x)+\alpha_{n+1}^l x_s$ is a weak* cluster point of
the sequence $(\sum_{i=1}^n\alpha_i^l (x_{k_i}-x)+\alpha_{n+1}^l(x_k -x))_{k=1}^\infty$ and for its norm we have
$$\begin{aligned}
\norm{\sum_{i=1}^n\alpha_i^l (x_{k_i}-x)+\alpha_{n+1}^l x_s}
& = \norm{\sum_{i=1}^n\alpha_i^l (x_{k_i}-x)}+\norm{\alpha_{n+1}^l x_s}
\\ & \geq(c-(1-2^{-n})\ep)\sum_{i=1}^n\betr{\alpha_i^l}+\betr{\alpha_{n+1}^l}(c-\frac\ep2)
\\ & >(c-(1-2^{-n})\ep)\sum_{i=1}^{n+1}\betr{\alpha_i^l}=c-(1-2^{-n})\ep.
\end{aligned}
$$
It follows that there is $k_{n+1}>k_n$ such that
$$\norm{\sum_{i=1}^{n+1}\alpha_i^l (x_{k_i}-x)}>c-(1-2^{-n})\ep$$ for all $l\in\{1,\dots,L\}$.
By a straightforward calculation using the choice of the $\alpha^l$ and the triangle inequality we get
that inequality \eqref{gl7}, with $n+1$ instead of $n$, holds for all $\alpha$ in the unit sphere of $\ell_1^{n+1}$
and hence for all elements of $\er^{n+1}$.

This finishes the construction. By Lemma~\ref{lemma-l1}(i) we get
$$\de{X}{x_{k_n}-x}\ge 2(c-\ep),$$
hence clearly
$$\de{X}{x_k}\ge\de{X}{x_{k_n}}=\de{X}{x_{k_n}-x}\ge 2(c-\ep).$$
As $\varepsilon\in(0,c)$ is arbitrary, we get \eqref{est-L1bis}.

\smallskip

Finally, let us prove \eqref{est-L1tres}: 
We take any subsequence $(x_{k_n})$ and observe that
$$2\operatorname{d}(\clu{X}{x_k},X)\leq2\dh(\clu{X}{x_{k_n}}, X)\leq\de{X}{x_{k_n}}$$
by \eqref{est-L1bis}. Then we can pass to the infimum over all $(x_{k_n})$.
This finishes the proof of the theorem.

\section{Proof of Example~\ref{exa1}}

For $n\in\en$ set $X_n=\ell_\infty^n$ and let $X$ be the $\ell_1$-sum of all the spaces $X_n$, $n\in\en$.
Then $X$ is $L$-embedded by \cite[Proposition IV.1.5]{hawewe}.

Further, let $e_1^n,\dots,e^n_n$ be the canonical basic vectors of $X_n$ and let $(x_k)$ be the sequence in $X$ containing subsequently these basic vectors, i.e., the sequence
$$ e_1^1,e_1^2,e_2^2,e_1^3,e_2^3,e_3^3,e_1^4,\dots,e_4^4,\dots$$
Then we have $\wde{x_k}=2$ as each subsequence of $(x_k)$ contains a further subsequence isometrically equivalent to the canonical basis of $\ell_1$.

It remains to show that $\operatorname{d}(\clu{X}{x_k},X)=0$. To do so, it is enough to prove that $0$ is a weak cluster point of the sequence $(x_k)$.
To verify this, we fix $g^1,\dots,g^m\in X^*$ and $\ep>0$. Let $K=\max\{\|g^1\|,\dots,\|g^m\|\}$.

The dual $X^*$ can be canonically identified with the $\ell_\infty$-sum of the spaces $X_n^*$, $n\in\en$. Moreover, $X_n^*$ is canonically isometric to $\ell_1^n$. Thus each $g\in X^*$ can be viewed as a bounded sequence $(g_n)_{n\in\en}$, where $g_n=(g_{n,j})_{j=1}^n\in\ell_1^n$ for each $n\in\en$.

We find $N\in\en$ such that $\frac K N<\ep$ and let $n\in\en$ be such that $n > m N$.
Let $k\in\{1,\dots,m\}$ be arbitrary. We have $\|g^k_n\|\le\|g^k\|\le K$. As
$\|g^k_n\|=\sum_{j=1}^n |g^k_{n,j}|$, the set
$$\{j\in \{1,\dots,n\}: |g^k_{n,j}|\geq \tfrac K N\}$$
has at most $N$ elements. It follows that the set 
$$\{j\in \{1,\dots,n\}: (\exists k\in\{1,\dots,m\}, |g^k_{n,j}|\geq \tfrac K N)\}$$
has at most $m N$ elements. As $n> m N$, there is some $j\in\{1,\dots,n\}$ such that $|g^k_{n,j}|< \tfrac K N<\ep$ for each $k\in\{1,\dots,m\}$.
It means that $|g^k(e_j^n)|<\ep$ for each $k\in\{1,\dots,m\}$. 

Since $e_j^n$ is an element of the sequence $(x_k)$, this completes the proof that $0$ is in the weak closure of the sequence, hence $0$ is a weak cluster point (as the sequence $(x_k)$ does not contain $0$).

\section{Proof of Example~\ref{exa2}}

We recall that $\beta\en$ is the \v{C}ech--Stone compactification of $\en$ and $M(\beta\en)$ is the space of all signed Radon measures on $\beta\en$ considered as the dual of $\ell_\infty$.

Let us fix $\alpha>0$ and consider the space 
\[
Y_\alpha=(\ell_1,\alpha\|\cdot\|_1)\oplus_1 (C[1,\omega],\|\cdot\|_\infty).
\]
Here $\|\cdot\|_1$ denotes the usual norm on $\ell_1$, $\omega$ is the first infinite ordinal, $C[1,\omega]$ stands for the space of all continuous functions on the ordinal interval $[1,\omega]$ and $\|\cdot\|_\infty$ is the standard supremum norm. 
Note that we have the following canonical identifications:
$$\begin{aligned}
Y_\alpha^*&=(\ell_\infty,\tfrac1\alpha\|\cdot\|_\infty)\oplus_\infty (\ell_1[1,\omega],\|\cdot\|_1), \mbox{ and}\\
Y_\alpha^{**}&=(M(\beta\en),\alpha\|\cdot\|_{M(\beta\en)})\oplus_1 (\ell_\infty[1,\omega],\|\cdot\|_\infty).
\end{aligned}$$

For $k\in\en$, let $x_k=(e_k,\chi_{[k,\omega]})\in Y_\alpha$, where $e_k$ denotes the $k$-th canonical basic vector in $\ell_1$ and $\chi_{[k,\omega]}$ is the characteristic function of the interval $[k,\omega]$. 
Let $X_\alpha$ be the closed linear span of the set $\{x_k:k\in\en\}$.
We observe that
\begin{equation}\label{Xalpha}X_\alpha=\left\{((\eta_k),f)\in Y_\alpha:
f(n)=\sum_{k=1}^n\eta_k\mbox{ for all }n\in\en\right\}.\end{equation}
Indeed, the set on the right-hand side is a closed linear subspace of $Y_\alpha$ containing $x_k$ for each $k\in\en$. This proves the inclusion `$\subset$'. To prove the converse one, let us take any point $((\eta_k),f)$ in the set on the right-hand side. Since $(\eta_k)\in\ell_1$, we get
$$((\eta_k),f)=\sum_{k=1}^\infty \eta_k x_k\in X_\alpha$$
as the series is absolutely convergent.

It follows that for each $((\eta_k),f)\in X_\alpha$ we have
$$\alpha\|(\eta_k)\|_1\le \|((\eta_k),f)\| \le (\alpha+1)\|(\eta_k)\|_1,$$
hence $X_\alpha$ is isomorphic to $\ell_1$. More precisely, the projection on the first coordinate is an isomorphism onto $\ell_1$. In particular,
$X_\alpha$ has the Schur property (and thus it is weakly sequentially complete).

We further observe that $X_\alpha^{**}$ is canonically identified with the weak* closure of $X_\alpha$ in $Y_\alpha^{**}$, thus
\begin{multline}\label{Xalpha**} X_\alpha^{**}=\{(\mu,f)\in M(\beta\en)\times \ell_\infty[1,\omega]: \\
(\forall n\in\en: f(n)=\mu\{1,\dots,n\})\mbox{ and } f(\omega)=\mu(\beta\en)\}.\end{multline}
Indeed, the set on the right-hand side is  a weak* closed linear subspace of $Y_\alpha^{**}$ containing $X_\alpha$, which proves the inclusion `$\subset$'.
To prove the converse one let us fix $(\mu,f)$ in the set on the right-hand side.
Take a bounded net $(u_\tau)$ in $\ell_1$ which weak* converges to $\mu$. For each $\tau$ there is a unique $f_\tau\in C[1,\omega]$ such that $(u_\tau,f_\tau)\in X_\alpha$. Then $(f_\tau)$ is clearly a bounded net in $\ell_\infty[1,\omega]$. Moreover, we will show that $(f_\tau)$ weak* converges to $f$. Since the weak* topology  on bounded sets coincides with the topology of pointwise convergence, it suffices to show that $f_\tau$ pointwise converge to $f$. Indeed,
$$\begin{aligned}
f_\tau(n)&=\sum_{k=1}^n u_\tau(k) \to \mu\{1,\dots,n\}=f(n),\mbox{ for each }n\in\en,\\
f_\tau(\omega)&= \sum_{k=1}^\infty u_\tau(k)\to\mu(\beta\en)=f(\omega).
\end{aligned}$$  
It follows that $X_\alpha$ is $1$-complemented in its bidual. To show that we set
$$P(\mu,f)=((\mu\{k\}),f-\mu(\beta\en\setminus\en)\cdot\chi_{\{\omega\}}),\qquad (\mu,f)\in X_\alpha^{**}.$$
Then $P$ is a projection of $X_\alpha^{**}$ onto $X_\alpha$ of norm one. Indeed, if $(\mu,f)\in X_\alpha$, then $\mu(\beta\en\setminus\en)=0$ and hence $P(\mu,f)=(\mu,f)$. Further, by \eqref{Xalpha} and \eqref{Xalpha**} we get that $P(\mu,f)\in X_\alpha$ for each $(\mu,f)\in X_\alpha^{**}$. Thus $P$ is a projection onto $X_\alpha$. To show it has norm one, it is enough to observe that, given $(\mu,f)\in X_\alpha^{**}$, we have $\|(\mu\{k\})\|_{\ell_1}\le\|\mu\|$, and that $f-\mu(\beta\en\setminus\en)\cdot\chi_{\{\omega\}}$ is a continuous function on $[1,\omega]$ coinciding on $[1,\omega)$ with $f$ and so $\|f-\mu(\beta\en\setminus\en)\cdot\chi_{\{\omega\}}\|_\infty\le\|f\|_\infty$.

Further, for the sequence $(x_k)$, its weak$^*$ cluster points in 
$X_\alpha^{**}$
are equal to
\[
\{(\ep_t, \chi_{\{\omega\}}): t\in \beta\en\setminus\en\},
\]
where $\ep_t$ denotes the Dirac measure at a point $t\in\beta\en$.

We claim that, for our sequence $(x_k)$, we have 
\begin{equation}
\label{est-alpha}
\dd(\clu{X_\alpha}{x_k}, X_\alpha)\geq \frac12\quad\text{and}\quad \de{X_\alpha}{x_k}=2\alpha.
\end{equation}
To see the first inequality, we use the fact that the distance of any weak$^*$ cluster point of $(x_k)$ from $X_\alpha$ is at least $\operatorname{d}(\chi_{\{\omega\}},C[1,\omega])=\frac12$. On the other hand, if $t,t'\in\beta\en\setminus\en$ are distinct, then 
\[
\|(\ep_t,\chi_{\{\omega\}})-(\ep_{t'},\chi_{\{\omega\}})\|_{X_\alpha^{**}}=\|(\ep_t-\ep_{t'},0)\|_{X_\alpha^{**}}
=\alpha\|\ep_t-\ep_{t'}\|_{M(\beta\en)}=2\alpha.
\]
This verifies \eqref{est-alpha}.

Now we use the described procedure to construct the desired space $X$.
For $n\in\en$, let $\alpha_n=\frac1n$ and let $X_{\frac1n}$ be the space constructed for $\alpha_n$. Let 
\[
X=\left(\sum_{n=1}^\infty X_{\frac1n}\right)_{\ell_1}
\]
be the $\ell_1$-sum of the spaces $X_{\frac1n}$. We claim that $X$ is the  required space. 

First, since each $X_{\frac1n}$ has the Schur property, $X$, as their $\ell_1$-sum, possesses this property as well (this follows
by a straightforward modification of the proof that $\ell_1$ has the Schur property, see \cite[Theorem~5.19]{fhhmpz}).
Hence $X$ is weakly sequentially complete. 

Further, observe that 
$$X^*= \left(\sum_{n=1}^\infty X_{\frac1n}^*\right)_{\ell_\infty}
\mbox{\qquad and\qquad}X^{**}\supset\left(\sum_{n=1}^\infty X_{\frac1n}^{**}\right)_{\ell_1}.$$
Note that the latter space is not equal to $X^{**}$ but it is  $1$-complemented in $X^{**}$ (cf.\ the proof of \cite[Proposition IV.1.5]{hawewe}).
Now it follows that $X$ is $1$-complemented in $X^{**}$.

Finally, fix $n\in\en$.  We consider a sequence $\wh{x}_k=(0,\dots,0,\stackrel{n\text{-th}}{x_k},0,\dots)$,
where the elements $x_k\in X_{\frac1n}$, $k\in\en$, are defined above.
Let $y=(0,\dots,0,\stackrel{n\text{-th}}{(\ep_t,\chi_{\{\omega\}})},0,\dots)$, where $t\in \beta\en\setminus \en$, be a weak$^*$ cluster point of
$(\wh{x}_k)$ in $X^{**}$.
Then, for any $z=(z(1),z(2),\dots)\in X$,
\[
\|y-z\|_{X^{**}}\geq \|(\ep_t,\chi_{\{\omega\}})-z(n)\|_{X_{\frac1n}^{**}}\geq \frac12
\]
by \eqref{est-alpha}. Hence 
\[
\dd(\clu{X}{\wh{x}_k}, X)\geq \frac12.
\]

On the other hand,
\[
\de{X}{\wh{x}_k}= \de{X_{\frac1n}}{x_k}=\frac2n,
\]
again by \eqref{est-alpha}. From this observation the conclusion follows.

\section{Final remarks}

Up to now we have tacitly assumed that we are dealing with real Banach spaces.
In fact, our proofs work for real spaces but all the results can be easily transferred to complex spaces as well. Let us indicate how to see this.

Let $X$ be a complex Banach space. Denote by $X_R$ the same space considered over the field of real numbers (i.e., we just forget multiplication by imaginary numbers). Let $\phi: X^*\to (X_R)^*$ be defined by 
$$\phi(x^*)(x)=\Re x^*(x),\qquad x^*\in X^*, x\in X.$$
It is well known that $\phi$ is a real-linear isometry of $X^*$ onto $(X_R)^*$.
Let us define a mapping $\psi: X^{**}\to (X_R)^{**}$ by the formula
$$\psi(x^{**})(y^*)=\Re x^{**}(\phi^{-1}(y^*)),\qquad x^{**}\in X^{**}, y^*\in
(X_R)^*.$$
It is easy to check that the mapping $\psi$ satisfies the following properties:
\begin{itemize}
	\item[(i)] $\psi$ is a real-linear isometry of $X^{**}$ onto $(X_R)^{**}$.
	\item[(ii)] $\psi$ is a weak*-to-weak* homeomorphism.
	\item[(iii)] $\psi(X)=X_R$. 
\end{itemize}
It follows that for any sequence in $X$ all the quantities in question (i.e.,
$\delta$, $\widetilde\delta$, $\operatorname{d}$ and $\dh$) are the same with respect to $X$ and with respect to $X_R$.
(Recall that $\delta$ is defined as the diameter of weak* cluster points, which has good sense in a complex space as well,
even though in the complex case only the second formula of \eqref{gl} works.)
If, moreover, we observe that $X_R$ is $L$-embedded whenever $X$ is $L$-embedded,
we conclude that Theorem~\ref{positive} is valid for complex spaces as well.

As for Examples~\ref{exa1} and~\ref{exa2}, it is clear that they work also in the complex setting -- we can just consider
complex versions of the respective spaces.

\medskip

We finish by recalling that G.\ Godefroy's question, for which Banach spaces \eqref{est-L1} holds, remains open. In particular,
the following question seems to be open.

\begin{question} Let $X$ be a Banach space which is a $u$-summand in its bidual,
i.e., there is a projection $P:X^{**}\to X$ with $\|I-2P\|=1$. Does \eqref{quantitative} hold for $X$ for some $C>0$?
\end{question}

We conjecture that the space from Example~\ref{exa2}, although it is $1$-complemented in its bidual, is not a $u$-summand. At least the  projection we have constructed does not work.

\end{document}